\documentclass{amsart}

\title{Proximally well-monotone covers and $QH$-singularity}
\author{Tom Vroegrijk}
\date{}

\usepackage{amsfonts}
\usepackage{amsmath}

\newtheorem{prop}{Proposition}
\newtheorem{cor}{Corollary}

\theoremstyle{definition}

\newtheorem{defn}{Definition}

\newcommand{\cf}{\operatorname{cf}}
\newcommand{\tb}{\operatorname{b}}

\begin{document}

\maketitle


\begin{abstract}
In this paper we use a certain class of well-monotone covers on a quasi-uniform space $(X, \mathcal{U})$ to investigate whether there are quasi-uniformities $\mathcal{V}$ that are distinct from $\mathcal{U}$, but have the property that the associated Hausdorff quasi-uniformities $\mathcal{U}_H$ and $\mathcal{V}_H$ on the hyperspace of $X$ have the same underlying topologies.
\end{abstract}

\section{Introduction}

The Hausdorff distance, a distance function on the collection of subsets of a metric space, was first introduced by Hausdorff in \cite{Hausdorff}. A slight variation of this concept, however, had been defined earlier by Pompeiu in \cite{Pompeiu}. The Hausdorff distance has had many applications in various branches of mathematics since its first appearance and, even though it has been around for quite a while, it is still an essential tool in pattern detection and face recognition software.

Similar to the way one can construct the Hausdorff distance on the hyperspace of a metric space $(X, d)$, one can define a quasi-uniformity $\mathcal{U}_H$ given a quasi-uniform space $(X, \mathcal{U})$. This filter of entourages $\mathcal{U}_H$ forms a quasi-uniformity on the collection $\mathcal{P}(X)$ of all subsets of $X$ and is referred to as the \emph{Hausdorff quasi-uniformity} or \emph{Hausdorff-Bourbaki quasi-uniformity} associated with $\mathcal{U}$.

Whenever two Hausdorff quasi-uniformities $\mathcal{U}_H$ and $\mathcal{V}_H$ induce the same topology the quasi-uniformities $\mathcal{U}$ and $\mathcal{V}$ are called $QH$-equivalent (see \cite{Cao}). When working in the symmetric setting (and thus only studying uniformities) this concept is also called $H$-equivalence. Smith \cite{Smith} and Ward \cite{Ward1, Ward2} laid the foundations for the analysis of $H$-equivalent uniformities. Among other results, Smith proved that metrisable uniformities and totally bounded uniformities are $H$-singular, i.e. there are no distinct uniformities that are $H$-equivalent to them. Distinct uniformities $\mathcal{U}$ and $\mathcal{V}$ on a set $X$ can generate Hausdorff quasi-uniformities with the same underlying topology on the hyperspace $\mathcal{P}(X)$. Nevertheless, Poljakov \cite{Poljakov} has shown that whenever $(\mathcal{U}_H)_H$ and $(\mathcal{V}_H)_H$ define the same topology on $\mathcal{P}(\mathcal{P}(X))$, the uniformities $\mathcal{U}$ and $\mathcal{V}$ must coincide.

In the asymmetric case the concept of $QH$-equivalence behaves somewhat differently. First of all Cao et al. \cite{Cao} gave examples of quasi-uniform spaces that were either quasi-metrisable or totally bounded, but not $QH$-singular. Furthermore, K\"unzi showed in \cite{Kunzi} that Poljakov's result does not hold in the asymmetric case. In this paper we will investigate $QH$-singularity by means of proximally well-monotone covers. In particular we will give some results on quasi-metrisable spaces and totally bounded spaces that are $QH$-singular.

\section{Preliminaries}

Let $X$ be a set and $U, V\subseteq X\times X$ relations on $X$. For an $x\in X$ we define $U(x)$ as $\{y\in X\, |\, (x, y)\in U\}$. The relation $V\circ U$ 
contains all $(x, z)$ for which there is a $y\in X$ such that $y\in U(x)$ and $z\in V(y)$. We will denote $U\circ U$ as $U^2$ and $U\circ U^{n}
$ as $U^{n+1}$ whenever $n\geq 2$. A sequence $(U_n)_n$ of relations on $X$ is called a \emph{normal sequence} iff $U^2_{n+1}\subseteq U_n$ for each $n$.

A filter $\mathcal{U}$ on $X\times X$ is called a \emph{quasi-uniformity} iff it has the following properties:
\begin{enumerate}
\item $\forall x\in X\, \forall U\in \mathcal{U}:\, (x,x)\in U$,
\item $\forall U\in\mathcal{U}\, \exists V\in\mathcal{U}: \, V^2\subseteq U$.
\end{enumerate}
The elements of a quasi-uniformity $\mathcal{U}$ will be called \emph{entourages}. The pair $(X, 
\mathcal{U})$ is a \emph{quasi-uniform space}. Let $U$ be a relation on a quasi-uniform space $(X, \mathcal{U})$. We 
will say that $U$ is \emph{normal} with respect to $\mathcal{U}$ iff there is a normal sequence $(U_n)_n$ of entourages in $\mathcal{U}$ where $U_0$ is equal to $U$. Whenever it is clear in which quasi-uniformity we are working we will simply say that $U$ is normal. It is clear that each entourage of a quasi-uniformity is normal.

Each quasi-uniformity $\mathcal{U}$ has an underlying topology $\tau(\mathcal{U})$. In this topology the 
neighbourhoodfilter of a point $x$ is generated by the sets $U(x)$ with $U\in \mathcal{U}$. The quasi-uniformity $\mathcal{U}^{-1}$ is called the 
\emph{conjugate  of $\mathcal{U}$} and consists of all entourages $U^{-1}$, where $$U^{-1}=\{(y, x)\,|\, (x, y)\in \mathcal{U}\}.$$ The filter $\mathcal{U}^s$ 
that is generated by the relations $U\cap U^{-1}$ is a quasi-uniformity that is called the \emph{symmetrisation of $\mathcal{U}$}. If $\mathcal{U}$ is equal to its own symmetrisation, then $\mathcal{U}$ is called a \emph{uniformity} and $(X, \mathcal{U})$ a \emph{uniform space}. For an extensive 
monograph on quasi-uniform spaces we refer the reader to \cite{Fletcher}.

The set of all subsets of $X$ will be denoted as $\mathcal{P}(X)$. For a subset $A\in \mathcal{P}(X)$ and an entourage $U\in \mathcal{U}$ we define $U(A)
$ as $$\bigcup_{x\in A} U(x).$$ A subset $A$ of a quasi-uniform space is \emph{uniformly isolated} iff there is an entourage $U$ with the property $U(A)$ equals $A$. A quasi-uniform space $(X, \mathcal{U})$ in which each uniformly isolated subset is either empty of the entire space $X$ will be called \emph{uniformly connected}. For a given quasi-uniform space $(X, \mathcal{U})$ the relation $\ll$ on $\mathcal{P}(X)$ is defined such that $A\ll B$ iff 
there is a $U\in\mathcal{U}$ with the property that $U(A)\subseteq B$.

For any relation $U$ on $X$ we define $$U_+=\{(A, B)\in \mathcal{P}(X)\times \mathcal{P}(X)\, | \, 
B\subseteq U(A)\}$$ and $$U_-=\{(A, B)\in \mathcal{P}(X)\times \mathcal{P}(X)\, | \, A\subseteq U^{-1}(B)\}.$$ If $(X, \mathcal{U})$ is a quasi-uniform space, then the filter generated by the sets $U_-$ is a 
quasi-uniformity on $\mathcal{P}(X)$ that we will call the lower Hausdorff quasi-uniformity. Analogously, the sets $U_+$ generate the upper Hausdorff 
quasi-uniformity on $\mathcal{P}(X)$. We will denote the intersection $U_-\cap U_+$ as $U_H$. The Hausdorff quasi-uniformity $\mathcal{U}_H$ on the 
hyperspace $\mathcal{P}(X)$ is the filter that is generated by the sets $U_H$.

If $\mathcal{U}$ and $\mathcal{V}$ are two quasi-uniformities on a set $X$, then we say that $\mathcal{V}$ is \emph{$QH$-finer} than $\mathcal{U}$ iff $
\tau(\mathcal{U}_H)\subseteq \tau(\mathcal{V}_H)$. If the topologies $\tau(\mathcal{U}_H)$ and $\tau(\mathcal{V}_H)$ are equal, then we say that $\mathcal{U}$ and 
$\mathcal{V}$ are \emph{$QH$-equivalent}. A quasi-uniformity $\mathcal{U}$ is called \emph{$QH$-singular} iff there is no quasi-uniformity, other than $
\mathcal{U}$ itself, that is $QH$-equivalent to $\mathcal{U}$.

Let $(X, q)$ be a quasi-pseudometric space. The conjugate quasi-metric $q^{-1}$ on $X$ 
is defined such that $q^{-1}(x, y)$ is equal to $q(y, x)$  for all $x, y\in X$. The metric $q\vee q^{-1}$ is generally denoted as $q^s$. Its underlying uniformity 
is equal to the symmetrisation of the underlying quasi-uniformity of $q$. A quasi-pseudometric space $(X, q)$ will be called $QH$-singular iff its underlying quasi-uniformity is $QH$-singular For a set $A\subseteq X$ and an $x\in X$ the value $q(x, A)$ is defined as $\inf\{q(x, a)| a\in A\}$.

\section{Proximally Well-monotone Covers}

In this section we will introduce proximally well-monotone covers of a quasi-uniform space. We will use these covers to construct quasi-uniformities that are 
$QH$-equivalent to a given quasi-uniformity.

The following results generalise a theorem of Albrecht \cite{Albrecht} on the comparison of Hausdorff uniformities. We will use them extensively throughout 
this text to investigate $QH$-equivalence of quasi-uniformities.

\begin{prop}\label{qhfiner}
If $U$ and $V$ are relations on $X$ and $A$ a subset of $X$, then $U_H(A)\subseteq V_H(A)$ iff $U(A)\subseteq V(A)$ and for each $x\in A$ there is an $y\in A$ with the property $U(y)
\subseteq V(x)$.
\begin{proof}
First we will prove the sufficiency of this proposition. Suppose $B$ is an element of $U_H(A)$, then we know that by definition $B\subseteq U(A)$ and thus $B\subseteq V(A)$. Now take an 
$x\in A$ and a $y\in A$ such that $U(y)\subseteq V(x)$. By assumption we can find a $z\in B$ with $y\in U^{-1}(z)$. This means that $z\in U(y)\subseteq 
V(x)$ and that $A\subseteq V^{-1}(B)$. We can conclude that $B$ is an element of $V_H(A)$.

Conversely, assume that $U_H(A)\subseteq V_H(A)$. It is clear that the conditions stated in the proposition are true if $A$ is empty. So let us assume that $A$ is non-empty and take an $x\in U(A)$. It is clear that the set $A\cup \{x\}$ is an element of $U_H(A)$ and thus of $V_H(A)$. This implies that we can find a $y\in A$ 
such that $x\in V(y)$ and therefore $x\in V(A)$. Hence we have that $U(A)\subseteq V(A)$. Now let $x$ be an element of $A$ and assume that there is no $y\in A$ for which $U(y)\subseteq V(x)$. 
Choose for each $y\in A$ a $z_y\in U(y)\setminus V(x)$ and define $B$ as $\{z_y|y\in A\}$. By construction we have that $B$ is an element of $U_H(A)$ 
and thus of $V_H(A)$. This, however, cannot be possible since there is no $z_y\in B$ for which $z_y\in V(x)$.
\end{proof}
\end{prop}

From this proposition we immediately obtain the following result that describes when a quasi-uniformity $\mathcal{U}$ is $QH$-finer than $\mathcal{V}$.

\begin{cor}
$\mathcal{U}$ is $QH$-finer than $\mathcal{V}$ iff for each $A\subseteq X$ and $V\in\mathcal{V}$ there is a $U\in\mathcal{U}$ such that $U(A)\subseteq 
V(A)$ and for each $x\in A$ there is an $y\in A$ with the property $U(y)\subseteq V(x)$.
\end{cor}

Let $\mathcal{G}$ be a cover of $X$ that is well-ordered for the inclusion order and that does not have a maximal element. From here on we will denote the 
successor of an element $G$ in such a collection for this well-order as $G^+$. Because we assumed that $\mathcal{G}$ does not have a maximal element 
we know that each $G$ in fact has a successor. From this point on we will assume that $(X, \mathcal{U})$ is a quasi-uniform space. 

\begin{defn}
We will call a family $(U_G)_{G\in\mathcal{G}}$ of entourages that is indexed by a collection $\mathcal{G}
$ of subsets of $X$ \emph{decreasing} if and only if $U_{G'}\supseteq U_G$ whenever $G'\subseteq G$. It will be called \emph{normally decreasing} iff we can find 
decreasing families $(U_{(G, n)})_{G\in\mathcal{G}}$ such that for each $G\in\mathcal{G}$ we have that $(U_{(G, n)})_n$ is a normal sequence for $U_G$.
\end{defn}

\begin{defn}
A cover $\mathcal{G}$ of $X$ will be called a \emph{proximally well-monotone cover} iff all of its elements are non-empty, it is well-ordered for the inclusion order, does not have a maximal 
element and there is a normally decreasing family $(U_G)_{G\in\mathcal{G}}$ of entourages such that $U_G(G)\subseteq G^+$ for each $G\in\mathcal{G}$. 
\end{defn}

A quasi-pseudometric $q$ will be called \emph{uniformly continuous} for a quasi-uniform space $(X, \mathcal{U})$ iff $\{(x, y)\in X\times X| q(x,y)<\epsilon\}$ is an entourage for each $\epsilon>0$. Each quasi-uniform space is uniquely defined by the collection of all uniformly continuous quasi-pseudometrics. The following result defines proximally well-monotone covers in terms of quasi-pseudometrics. For a set $A
\subseteq X$ and an $\epsilon>0$ we define $A^{\epsilon}_q$ as $\{x\in X| \exists a\in A: q(a, x)<\epsilon\}$.

\begin{prop}
A cover $\mathcal{G}$ of $X$ is proximally well-monotone iff it all of its elements are non-empty, is well-ordered for the inclusion order, does not have a maximal element and there is a 
family $(q_G)_{G\in \mathcal{G}}$ of uniformly continuous quasi-pseudometrics and an $\epsilon>0$ such that $q_G\leq q_{G'}$ whenever $G\subseteq G'$ and $G^{\epsilon}
_{q_G}\subseteq G^+$ for each $G\in \mathcal{G}$.
\begin{proof}
Suppose that $\mathcal{G}$ is a cover that satisfies the conditions stated in the proposition. Define $U_{(G,n)}$ as $\{(x,y)\in X\times X\, |\, q_G(x, y)<2^{-n}
\epsilon\}$. Since each $q_G$ is uniformly continuous we have that each $U_{(G, n)}$ is an entourage. It is clear from the definition that for each $n$ the family $(U_{(G, n)})_{G\in\mathcal{G}}$ is decreasing and that each $(U_{(G,n)})_{n}$ is a normal sequence. 
Furthermore, we also have that for each $G\in\mathcal{G}$ the set $U_{(G, 0)}(G)$ is equal to $G_{q_G}^{\epsilon}$ and thus a subset of $G^+$. We can 
conclude that $\mathcal{G}$ is a proximally well-monotone cover.

Let $\mathcal{G}$ be a proximally well-monotone cover. By definition we can find for each $n$ a decreasing family $(U_{(G, n)})_{G\in\mathcal{G}}$ such 
that for each $G$ we have that $U_{(G, 0)}(G)\subseteq G^+$ and $(U_{(G, n)})_{n}$ is a normal sequence. In \cite{Kelley} Kelley showed that given a sequence of entourages $(V_n)_n$ of a quasi-uniform space with the property that $V^4_{n+1}\subseteq V_n$ for each $n$, we can always construct a quasi-pseudometric $q$ such that $$V_{n+1}\subseteq \{(x, y)\in X\times X\, |\, q(x, y)< 
2^{-n}\}\subseteq V_{n}.$$ Since for each $G\in \mathcal{G}$ the sequence $(U_{(G, 2n)})_n$ satisfies this property this means that we can construct a quasi-pseudometric $q_G$ that satisfies $$U_{(G, 2n+2)}\subseteq \{(x, y)\in X\times X\, |\, q_G(x, y)< 
2^{-n}\}\subseteq U_{(G, 2n)}$$ for each $n$. Take $G\subseteq G'\in \mathcal{G}$. According to Kelley's construction, the quasi-pseudometric $q_G$ is the least upper bound of all quasi-pseudometrics $q$ for which $q(x, y)$ is less than or equal to $$\inf \{2^{-n} | (x,y)\in U_{(G, 2n)}\}.$$ Suppose $G\subseteq G'$. For every $n$, since the family $(U_{(G, n)})_{G\in \mathcal{G}}$ is decreasing, we have that $U_{(G, n)}\supseteq U_{(G', n)}$. Hence, for every $x, y \in X$, we get $$\inf \{2^{-n} | (x,y)\in U_{(G, 2n)}\}\leq \inf \{2^{-n} | (x,y)\in U_{(G', 2n)}\}$$ and therefore $q_G\leq q_{G'}$. By definition we have that $G^1_{q_G}\subseteq U_{(G, 0)}(G)\subseteq G^+$ for each $G\in \mathcal{G}$. This means that the conditions 
stated in the proposition are satisfied.
\end{proof}
\end{prop}

\begin{prop}\label{countable}
If $(G_n)_n$ is a strictly increasing sequence of non-empty sets that cover $X$ with the property that $G_n\ll G_{n+1}$ for each $n$, then the collection of all $G_n$ is 
a proximally well-monotone cover of $X$.
\begin{proof}
That $(G_n)_n$ is well-ordered for the inclusion order and does not have a maximal element follows from our assumptions. Choose a sequence $(U_n)_n
$ of entourages such that $U_n(G_n)\subseteq G_{n+1}$. Without loss of generality we can assume that $(U_n)_n$ is decreasing. Let $(V_{(n, m)})_m$ be 
a normal sequence for $U_n$. Define $W_{(n, m)}$ as $V_{(0, m)}\cap V_{(1, m)}\ldots \cap V_{(n, m)}.$ It is clear that each $(W_{(n, m)})_m$ is a normal 
sequence for $U_n$ and that the sequence $(W_{(n, m)})_n$ is decreasing. This means that $(U_n)_n$ is normally decreasing and that the collection of 
all sets $G_n$ is a proximally well-monotone cover.
\end{proof}
\end{prop}

\begin{defn}
Let $\mathcal{G}$ be a proximally well-monotone cover of $X$. Since each $x\in X$ is contained in one of the elements of $\mathcal{G}$, so we can define the relation 
$U_{\mathcal{G}}$ on $X$ such that $U_{\mathcal{G}}(x)$ is equal to $G^+$, where $G$ is the smallest element of $\mathcal{G}$ that contains $x$.
\end{defn}

\begin{prop}\label{normalentourage}
Let $(X, \mathcal{U})$ be a uniformly connected quasi-uniform space. If $\mathcal{G}$ is a proximally well-monotone cover, then we can find a proximally well-monotone cover $\mathcal{G}^*$ such that $(U_{\mathcal{G}
^*})^2\subseteq U_{\mathcal{G}}$.
\begin{proof}
Since $\mathcal{G}$ is proximally well-monotone we can find a normally decreasing family $(U_G)_{G\in\mathcal{G}}$ such that $(U_G)^2(G)
\subseteq G^+$. Define $\mathcal{G}^*$ as the union of $\mathcal{G}$ and $\{U_G(G)|G\in\mathcal{G}\}$. It is clear from the construction that this set only contains non-empty subsets, cannot have a maximal element and that it is still a cover of $X$.

The set $\mathcal{G}^*$ is well-ordered for the inclusion order. If $\mathcal{H}$ is a non-empty subset of $\mathcal{G}^*$ the set $\{G\in\mathcal{G}|G\in\mathcal{H}
\textrm{ or }U_G(G)\in\mathcal{H}\}$ has a minimum $G_{\mathcal{H}}$ in $\mathcal{G}$. Suppose $G_{\mathcal{H}}\in\mathcal{H}$. Any other element in $\mathcal{H}$ is equal to 
either $G$ or $U_G(G)$ for some $G\in\mathcal{G}$. Since $G_{\mathcal{H}}\subseteq G\subseteq U_G(G)$ we find that $G_{\mathcal{H}}$ is the minimum of the set $\mathcal{H}
$. Let us now assume that $G_{\mathcal{H}}\not\in \mathcal{H}$. This yields that $U_{G_{\mathcal{H}}}(G_{\mathcal{H}})\in \mathcal{H}$. We know that for each $G\in \mathcal{H}$ that is also 
contained in $\mathcal{G}$ we have that $G_{\mathcal{H}}\subseteq G$. Since we assumed that $G_{\mathcal{H}}$ is not contained in $\mathcal{H}$ we know that $G$ must be 
strictly larger than $G_{\mathcal{H}}$. Because $\mathcal{G}$ is well-ordered this yields that $(G_{\mathcal{H}})^+\subseteq G$ and thus we have $U_{G_{\mathcal{H}}}(G_{\mathcal{H}})\subseteq (G_{\mathcal{H}})^+
\subseteq G$. On the other hand if $U_G(G)\in \mathcal{H}$ with $G\in\mathcal{G}$, then $G_{\mathcal{H}}\subseteq G$. This implies that $U_{G_{\mathcal{H}}}(G_{\mathcal{H}})\subseteq 
U_{G}(G)$. In case $G$ equals $G_{\mathcal{H}}$ this is trivially true. If $G$ is strictly larger than $G_{\mathcal{H}}$ we find $U_{G_{\mathcal{H}}}(G_{\mathcal{H}})\subseteq (G_{\mathcal{H}})^+\subseteq G \subseteq 
U_G(G)$. We can conclude that $U_{G_{\mathcal{H}}}(G_{\mathcal{H}})$ is the minimum of $\mathcal{H}$.

So far we have established that $\mathcal{G}^*$ is a well-ordered family of non-empty sets without maximal element. To prove that it is also proximally well-monotone we take an arbitrary $G^*\in \mathcal{G}^*$. If $G^*$ is not contained in $\mathcal{G}$, then by definition this means that $G^*$ is equal to $U_G(G)$ for some 
unique $G\in\mathcal{G}$. Define $U_{G^*}$ as $U_G$. If $G^*$ is an element of $\mathcal{G}$, then $U_{G^*}$ is just the entourage that we defined earlier. It is clear that $(U_{G^*})_{G^*\in \mathcal{G^*}}$ is a normally decreasing family. If $G^*$ is equal to $U_G(G)$ for some $G\in\mathcal{G}$, then the 
successor of $G^*$ in $\mathcal{G}^*$ must contain $G^+$. Suppose that it is not equal to $G^+$. This implies that $U_G(G)=G^+$ and thus $U_G(U_G(G))=G^+$, but this means that $G^+$ is uniformly isolated and we assumed that $(X, \mathcal{U})$ is uniformly connected. By definition we have that $U_{G^*}(G^*)=U_G(U_G(G))\subseteq G^+$. On the other hand, if $G^*$ is an element in $\mathcal{G}$, then again from the uniform connectedness of $(X, \mathcal{U})$ it follows that its successor in $
\mathcal{G}^*$ is equal to $U_{G^*}(G^*)$. This means that $\mathcal{G}^*$ is a proximally 
well-monotone cover of $X$.

The only thing that is left to prove is the fact that $(U_{\mathcal{G}
^*})^2\subseteq U_{\mathcal{G}}$. Take an $x\in X$ and let $G_x$ be the smallest element in $\mathcal{G}^*$ that contains $x$. If $G_x$ is an element of $\mathcal{G}$, then $U_{\mathcal{G}
^*}(x)$ is equal to $U_{G_x}(G_x)$ and for a $y\in U_{\mathcal{G}^*}(x)$ we have that $U_{\mathcal{G}^*}(y)$ is a subset of the successor of $G_x$ in $\mathcal{G}
$. Since $U_{\mathcal{G}}(x)$ is equal to the successor of $G_x$ in $\mathcal{G}$ we find that $(U_{\mathcal{G}^*})^2(x)\subseteq U_{\mathcal{G}}(x)$. If $G_x$ is equal to $U_G(G)$ 
for some $G\in\mathcal{G}$, then $U_{\mathcal{G}^*}(x)$ is equal to the successor $G^+$ of $G$ in $\mathcal{G}$ and $U_{\mathcal{G}}(x)$ equals the 
successor of $G^+$ in $\mathcal{G}$. This means that for any $y\in U_{\mathcal{G}^*}(x)$ the set $U_{\mathcal{G}^*}(y)$ is a subset of $U_{G^+}(G^+)$ and 
thus $(U_{\mathcal{G}^*})^2(x)\subseteq U_{G^+}(G^+) \subseteq U_{\mathcal{G}}(x)$.
\end{proof}
\end{prop}

\begin{cor}
If $\mathcal{G}$ is a proximally well-monotone cover, then we can find a sequence $(\mathcal{G}_n)_n$ of proximally well-monotone covers such that $(U_{\mathcal{G}_n})_n$ is a normal sequence for $U_{\mathcal{G}}$.
\end{cor}

\begin{prop}\label{monotonecover}
Let $(X, \mathcal{U})$ be a quasi-uniform space. If $\mathcal{G}$ is a proximally well-monotone cover, then $(U_{\mathcal{G}})_H$ is an element of $
\mathcal{U}_H$.
\begin{proof}
Let $(U_G)_{G\in\mathcal{G}}$ be a normally decreasing sequence for $\mathcal{G}$ and take an $A\subseteq X$. We will verify that $U_{\mathcal{G}}$ 
satisfies the conditions given in proposition \ref{qhfiner}. Assume that $A$ is a subset of an element of $\mathcal{G}$. Define $G_1$ as the smallest 
element of $\mathcal{G}$ that intersects with $A$ and $G_2$ as the smallest set in $\mathcal{G}$ that contains $A$. It is clear that $G_1$ must be a subset 
of $G_2$. We now have $U_{G_2}(A)\subseteq U_{G_2}(G_2) \subseteq G_2^+ =U_{\mathcal{G}}(A)$. Take an $x\in A$. We know that $U_{G_2}\subseteq U_{G_1}$. For an 
arbitrary $y\in A\cap G_1$ we find that 
\begin{eqnarray*}
U_{G_2}(y) & \subseteq & U_{G_1}(y) \\
 & \subseteq & G_1^+ \\
 & \subseteq & G_2^+ \\
 & \subseteq & U_{\mathcal{G}}(x).
\end{eqnarray*}
We find that $U_{G_2}$ satisfies the conditions stated in proposition \ref{qhfiner} and thus we have that $(U_{G_2})_H(A)\subseteq U_{\mathcal{G}}(A)$.

Now assume that $A$ is not contained in any element of $\mathcal{G}$. Define $G_1$ again as the smallest element of $\mathcal{G}$ that intersects with 
$A$. For each $G\in \mathcal{G}$ we can find an $x\in A$ such that $x\not\in G$ and thus $G\subseteq U_{\mathcal{G}}(x)$. Because $\mathcal{G}$ is a 
cover of $X$ we obtain that $U_{\mathcal{G}}(A)$ is equal to $X$. This implies that $U_{G_1}(A)\subseteq U_{\mathcal{G}}(A)$. With the same arguments 
as we used above we can argue that for each $x\in A$ there is a $y\in A$ with $U_{G_1}(y)\subseteq U_{\mathcal{G}}(x)$ and therefore we can conclude that $(U_{G_1})_H(A)\subseteq U_{\mathcal{G}}(A)$.
\end{proof}
\end{prop}

\begin{prop}\label{monotonehaus}
Let $(X, \mathcal{U})$ be a quasi-uniform space. Let $\mathcal{G}$ be a proximally well-monotone cover with the property that for each $U\in \mathcal{U}$ 
and each $A\subseteq X$ that is not contained in any element of $\mathcal{G}$ there is a $G\in\mathcal{G}$ such that for each $x\in A\setminus G$ we can find:
\begin{itemize}
\item an $a'$ in $A\setminus G'$ with $x\in U(a')$ whenever $G'\in \mathcal{G}$,
\item an $a\in A\cap G$ such that $a\in U(x)$.
\end{itemize}
For each $V\in\mathcal{U}$ we have that $(U_{\mathcal{G}}\cap V)_H\in\mathcal{U}_H$.
\begin{proof}
Let $(U_G)_{G\in\mathcal{G}}$ be a normally decreasing family, take $V\in\mathcal{U}$ and $A\subseteq X$. First let us assume that $A$ is contained in an element $G$ of $\mathcal{G}$. By definition this means that $U_G(a)\subseteq U_{\mathcal{G}}(a)$ for each $a\in A$ and thus $(U_G\cap V)(a)\subseteq (U_{\mathcal{G}}\cap V)(a)$. Clearly this implies that $(U_G\cap V)_H(A)\subseteq (U_{\mathcal{G}}\cap V)_H(A)$.

From here on we will assume that $A$ is not contained in any element of $\mathcal{G}$. Choose $U\in\mathcal{U}$ such that 
$U^2\subseteq V$. We can assume that each $U_G$ is contained in $U$. If this is not the case we choose a normal sequence $(U_n)_n$ of entourages for $U$ and decreasing families $(U_{(G, n)})_{G\in \mathcal{G}}$ such that each $(U_{(G, n)})_n$ is a normal sequence for $U_G$. It is clear that for each $n$ the family $(U_n\cap U_{(G, n)})_{G\in \mathcal{G}}$ is again decreasing and that $(U_n\cap U_{(G, n)})_n$ is a normal sequence for $U\cap U_G$ for all $G\in \mathcal{G}$.

For the entourage $U$ that we chose we can find a $G\in\mathcal{G}$ with the 
properties that are stated above. We will prove that $(U_{\mathcal{G}}\cap V)_H(A)$ is a neighbourhood of $A$ for $\tau(\mathcal{V}_H)$ by using proposition \ref{qhfiner}. We will start with proving that $U_G(A)\subseteq (U_{\mathcal{G}}\cap V)(A)$. Let $y$ be an element of $U_G(A)$. Take $x\in A$ such that $y\in U_G(x)$. 
Let us suppose that $x$ is not contained in $G$. Choose a $G'\in\mathcal{G}$ that is larger than $G$ and contains both the elements $x$ and $y$. Such a 
$G'$ exists since $\mathcal{G}$ is a well-ordered cover of $X$. By assumption we can now find an $a'\in A\setminus G'$ such that $x\in U(a')$. Since $a'$ 
is not contained in $G'$ and $y$ is, we automatically have $y\in U_{\mathcal{G}}(a')$. Because $x\in U(a')$ and $y\in U_G(x)\subseteq U(x)$ we get $y\in V(a')$ and thus $y\in (U_{\mathcal{G}}\cap V)_H(a')$. In case $x$ is an element of $G$ we know that $U_G(x)\subseteq U_{G_x}(x)= U_{\mathcal{G}}(x)$, where $G_x$ is the smallest set in $\mathcal{G}$ that contains $x$ and thus $y\in U_{\mathcal{G}}(x)$. Since $U_G\subseteq V$ and $y\in U_G(x)$ this implies that $y$ is contained in $(U_{\mathcal{G}}\cap V)(A)$.

Take an $x\in A$. We need to prove that there is an element $a\in A$ such that $U_G(a)\subseteq (U_{\mathcal{G}}\cap V)(x)$. If $x$ is not contained 
in $G$, then we already saw that $U_G(x)\subseteq (U_{\mathcal{G}}\cap V)(x)$. This means that in this case we can choose $a$ to be equal to $x$. 
Suppose $x$ is not contained in $G$. By assumption we can find for such an $x$ an $a\in A\cap G$ with the property $a\in U(x)$. If we take a $z\in U_G(a)
$, then $z$ is also contained in $U_G(G)$ and therefore in $G^+$. Because $x\not\in G$ this means that each element of $\mathcal{G}$ that contains $x$ 
contains $G^+$ and thus contains $z$. This implies that $z\in U_{\mathcal{G}}(x)$. Moreover, it follows from $z\in U_G(a)$, $a\in U(x)$ and $U_G\subseteq U$ that $z\in V(x)$. Hence $U_G(a)\subseteq (U_{\mathcal{G}}\cap V)(x)$.
\end{proof}
\end{prop}

\begin{defn}
Let $(X, \mathcal{U})$ be a quasi-uniform space and $U$ an entourage. We will say that a set $A\subseteq X$ is \emph{$U$-small} iff $A\times A\subseteq U$. Let $
\alpha$ be a cardinal number. A set $A$ will be called $U$-$\alpha$-bounded iff $A$ is equal to the union of a family $(A_i)_{i\in I}$ of $U$-small sets such 
that $|I|<\alpha$. An \emph{$\alpha$-bounded set} will be a set that is $U$-$\alpha$-bounded for each entourage $U\in\mathcal{U}$. It is clear from these 
definitions that a set is totally bounded iff it is $\aleph_0$-bounded.
\end{defn}

Let $\mathcal{G}$ be a collection of subsets of $X$. Recall that a subset $\mathcal{G}'$ of $\mathcal{G}$ is cofinal iff each $G\in\mathcal{G}$ is contained 
in an element of $\mathcal{G}'$. The minimal cardinality of all cofinal subsets of $\mathcal{G}$ is called its cofinality and is generally denoted as $
\cf(\mathcal{G})$.

\begin{prop}\label{boundedhaus}
If $\mathcal{G}$ is a proximally well-monotone cover of $X$ with the property that for each $U\in \mathcal{U}$ there is a $G\in\mathcal{G}$ such that $X
\setminus G$ is $U$-$\cf(\mathcal{G})$-bounded, then $(U_{\mathcal{G}}\cap V)_H\in\mathcal{U}_H$ for each $V\in\mathcal{U}$.
\begin{proof}
We will prove that $\mathcal{G}$ satisfies the conditions stated in proposition \ref{monotonehaus}. Let $A$ be a subset of $X$ that is not contained in any element of $\mathcal{G}$ and take an entourage 
$U\in\mathcal{U}$. Choose $G_1\in\mathcal{G}$ such that $X\setminus G_1$ is $U$-$\cf(\mathcal{G})$-bounded. We can now write $A\setminus G_1$ as 
the union of a collection $(A_i)_{i\in I}$ of $U$-small sets with $|I|<\cf(\mathcal{G})$. Let $J$ be the set that contains all $i\in I$ for which $A_i$ is contained 
in an element of $\mathcal{G}$. Choose for each $i\in J$ a $G_i\in\mathcal{G}$ such that $A_i\subseteq G_i$. Since $|J|<\cf(\mathcal{G})$ there must be a 
$G_2\in \mathcal{G}$ that contains each $G_i$ with $i\in J$. If this were not the case, then $(G_i)_{i\in J}$ would be a cofinal collection with cardinality 
strictly smaller than $\cf(\mathcal{G})$. Without loss of generality we can assume that $G_2$ is larger than $G_1$.

Since for each $i\in I\setminus J$ we have that $A_i\not\subseteq G_2$ we can choose an $a_i$ in $A_i\setminus G_2$. There must be a $G_0\in
\mathcal{G}$ that contains all $a_i$ with $i\in I\setminus J$. Suppose this were not true, then we could choose a $G_i\in\mathcal{G}$ such that $a_i\in G_i$ 
for each $i\in I\setminus J$ and because no element $G$ of $\mathcal{G}$ contains each $a_i$ this collection $(G_i)_{i\in I\setminus J}$ would be cofinal. This is impossible because $|I\setminus J|<\cf(\mathcal{G})$. Since $G_0$ contains the elements $a_i$, that are not contained $G_2$, we know that it must contain $G_2$.

We will now show that this set $G_0$ satisfies the conditions that are stated in proposition \ref{monotonehaus}. Let $x$ be an element of $A\setminus G_0$. 
By construction $x$ must be an element of an $A_i$ with $i\in I\setminus J$. We know that $a_i$ is contained in $A\cap G_0$, but not in $G_1$. From the fact 
that $X\setminus G_1$ is $U$-small we can conclude that $a_i\in U(x)$. On the other hand, if $G'$ is an element of $\mathcal{G}$, then $A_i$ cannot be 
contained in $G'$ since $i\not\in J$. Choose an $a'\in A_i\setminus G'$. Both $x$ and $a'$ are elements of the complement of $G_1$, so $x\in U(a')$.
\end{proof}
\end{prop}

\section{$QH$-singularity}

We shall start this section with a result that, given a quasi-uniform space $(X, \mathcal{U})$, allows us to construct a quasi-uniformity that is strictly larger than, but still $QH$-equivalent to $\mathcal{U}$. We will use this result to discuss the cofinallity of proximally well-monotone covers on $QH$-singular spaces.

\begin{prop}\label{singular}
Let $(X, \mathcal{U})$ be a uniformly connected quasi-uniform space. If there exists a proximally well-monotone cover $\mathcal{G}$ of $X$ with the property that for each $U\in \mathcal{U}$ there is a $G\in\mathcal{G}$ such 
that $X\setminus G$ is $U$-$\cf(\mathcal{G})$-totally bounded, then $(X,\mathcal{U})$ is not $QH$-singular.
\begin{proof}
Suppose that there is a proximally well-monotone cover $\mathcal{G}$ of $X$ with the property that is stated in the proposition. We saw earlier on that we can find a sequence $(\mathcal{G}_n)_n$ of 
proximally well-monotone covers such that $(U_{\mathcal{G}_n})_n$ is a normal sequence for $U_{\mathcal{G}}$. Looking at the construction of this 
normal sequence (see proposition \ref{normalentourage}) we see that each $\mathcal{G}_n$ contains $\mathcal{G}$. This means that each $\mathcal{G}_n$ 
satisfies the conditions stated in proposition \ref{boundedhaus}. Let $\mathcal{U}^*$ be the quasi-uniformity generated by the entourages $U_{\mathcal{G}
_n}\cap V$ with $V\in \mathcal{U}$.

From proposition \ref{boundedhaus} we obtain that each of the entourages $(U_{\mathcal{G}_n}\cap V)_H$ is an element of $\mathcal{U}_H$. Because $\mathcal{U}_H$ is clearly a 
subset of $(\mathcal{U}^*)_H$ we can conclude that both quasi-uniformities are in fact equal. This means that $\mathcal{U}$ and $\mathcal{U}^*$ are $QH
$-equivalent.

We will now prove that $\mathcal{U}$ and $\mathcal{U}^*$ are distinct unformities by showing that $U_{\mathcal{G}}$ cannot be an element of $
\mathcal{U}$. Let us assume that $U_{\mathcal{G}}$ is in fact an entourage in the quasi-uniformity $\mathcal{U}$. If this is the case, then we can find a $V
\in\mathcal{U}$ such that $V(x)\subseteq U_{\mathcal{G}}(x)$ for each $x\in X$. This implies that $V(G)\subseteq U_{\mathcal{G}}(G)\subseteq G^+$ for 
each $G\in\mathcal{G}$. Choose $G\in\mathcal{G}$ such that $X\setminus G$ is $V$-$\cf(\mathcal{G})$-totally bounded. This means that we can write the 
complement of $G$ as the union of a collection $(A_i)_{i\in I}$ where each $A_i$ is $V$-small and $|I|<\cf(\mathcal{G})$. Choose an $a_i\in A_i$ 
for each $i\in I$ and let $G_i$ be an element of $\mathcal{G}$ with the property $a_i\in G_i$. Because $(G_i)_{i\in I}$ cannot be cofinal we know that there is a $G'$ that is larger than each $G_i$. Since each $G_i$ is larger than $G$ the same must be true for $G'$. Take an arbitrary element $x$ in the complement of $G'$. 
There is an $i\in I$ such that $x\in A_i$. This yields that $x\in V(a_i)$ and thus that $x\in V(G')\subseteq (G')^+$. This means that each element in the complement of $G'$ is contained in $(G')^+$. Since $G'$ is by definition a subset of $(G')^+$ we obtain that $(G')^+$ is equal to $X$. This is impossible since we assumed a proximally well-monotone family to not have a maximal element. Hence we can say that $U_{\mathcal{G}}$ is no element of $\mathcal{U}$ and that $\mathcal{U}$ is not $QH$-singular.
\end{proof}
\end{prop}

For a quasi-uniform space $(X,\mathcal{U})$ we will define $\tb(X,\mathcal{U})$ as the minimal cardinal $\alpha$ for which $(X, \mathcal{U})$ is $\alpha$-
bounded. The first result that we can derive from proposition \ref{singular} is that this cardinal number is in fact a strict upper bound for the cofinality of proximally well-monotone covers on a QH-singular quasi-uniform space $(X, \mathcal{U})$.

\begin{cor}\label{boundedness}
Let $(X,\mathcal{U})$ be a uniformly connected, $QH$-singular quasi-uniform space. If $\mathcal{G}$ is a proximally well-monotone cover, then $\cf(\mathcal{G})<\tb(X,
\mathcal{U})$.
\begin{proof}
Suppose that $\mathcal{G}$ is a proximally well-monotone cover with the property that $\cf(\mathcal{G})$ is at least $\tb(X,\mathcal{U})$. Since $X
\setminus G$ is $\cf(\mathcal{G})$-bounded for each $G\in \mathcal{G}$ we obtain that $\mathcal{G}$ satisfies the properties stated in proposition 
\ref{singular}. This would imply that $(X, \mathcal{U})$ is not $QH$-singular.
\end{proof}
\end{cor}

\begin{cor}\label{countablesingular}
Let $(X,\mathcal{U})$ be a uniformly connected, totally bounded quasi-uniform space. If $(X, \mathcal{U})$ is $QH$-singular and we have a $\ll$-increasing sequence $(G_n)_n
$ of sets that covers $X$, then there is some $G_n$ that is equal to $X$.
\begin{proof}
If no $G_n$ is equal to $X$, then the collection of all sets $G_n$ forms a proximally well-monotone cover by proposition \ref{countable}. Because both the cofinality of this cover and $\tb(X, \mathcal{U})$ are equal to $\aleph_0$, it would follow from corollary \ref{boundedness} that $(X, \mathcal{U})$ is not $QH$-singular.
\end{proof}
\end{cor}


We will now prove that a quasi-uniform space is $QH$-singular whenever its symmetrisation is compact.

\begin{prop}\label{uniformcoreflection}
If $\mathcal{U}_H=\mathcal{V}_H$, then $\tau(\mathcal{U}^s)=\tau(\mathcal{V}^s)$.
\begin{proof}
If $\mathcal{U}$ and $\mathcal{V}$ are $QH$-equivalent, then $\tau(\mathcal{U})$ is equal to $\tau(\mathcal{V})$ and $\tau(\mathcal{U}^{-1})$ is equal to $
\tau(\mathcal{V}^{-1})$ (see \cite[Remark 2]{Kunzi}). Because the topologies $\tau(\mathcal{U}^s)$ and $\tau(\mathcal{V}^s)$ are respectively equal to $
\tau(\mathcal{U})\vee \tau(\mathcal{U}^{-1})$ and $\tau(\mathcal{V})\vee \tau(\mathcal{V}^{-1})$ we obtain that $\tau(\mathcal{U}^s)$ and $\tau(\mathcal{V}
^s)$ coincide.
\end{proof}
\end{prop}

\begin{prop}
If $\tau(\mathcal{U}^s)$ is compact and $\mathcal{U}_H=\mathcal{V}_H$, then $\mathcal{V}\subseteq \mathcal{U}$.
\begin{proof}
Take a $V\in\mathcal{V}$ and choose $W\in \mathcal{V}$ such that $W^2\subseteq V$. Because $\mathcal{U}_H=\mathcal{V}_H$ we know that $
\tau(\mathcal{U})=\tau(\mathcal{V})$. This means that we can take for each $x\in X$ a $U_x\in\mathcal{U}$ such that $U_x^2(x)\subseteq W(x)$. Because $W^{-1}(x)$ is a neighbourhood of $x$ in $\tau(\mathcal{V}^{-1})$ and thus in $\tau(\mathcal{U}^{-1})$ we know that $(W^{-1}\cap U_x)(x)$ is a $\tau(\mathcal{U}^s)$-neighbourhood of $x$. Since $\tau(\mathcal{U}^s)$ is compact we can find $x_0, \ldots, x_n\in X$ such that $$X=\bigcup_{k=0}^n (W^{-1}\cap U_{x_k})(x_k).$$ Define $U$ as $U_{x_0}\cap \ldots \cap U_{x_n}$. Let $x$ be an arbitrary element of $X$ and take a $z\in U(x)$. Assume that $x$ is an element of $(W^{-1}\cap U_{x_k})(x_k)$ for some $0\leq k \leq n$. This implies that $z$ is an element of $U(U_{x_k}(x_k))$ and thus of $W(x_k)$. Because $x_k\in W(x)$ we have that $z\in W^2(x)\subseteq V(x)$. This means that $U(x)\subseteq V(x)$ and, because $x$ was arbitrary, that $U\subseteq V$.
\end{proof}
\end{prop}

\begin{cor}\label{symmcompact}
If $\tau(\mathcal{U}^s)$ is compact, then $\mathcal{U}$ is $QH$-singular.
\begin{proof}
Suppose that there is a quasi-uniformity $\mathcal{V}$ such that $\mathcal{U}_H=\mathcal{V}_H$. From the previous proposition we obtain that $\mathcal{V}
\subseteq \mathcal{U}$. Proposition \ref{uniformcoreflection} yields that $\tau(\mathcal{V}^s)$ is compact. Using the previous proposition once more, we 
get that $\mathcal{U}\subseteq \mathcal{V}$.
\end{proof}
\end{cor}

Using the results that we have obtained so far we can give a complete characterisation of $QH$-singularity for totally bounded, quasi-metric spaces.

\begin{prop}\label{dense}
Let $(X, q)$ be a uniformly connected quasi-metric space and $Y$ a non-empty subset of $X$. If there is an $x\in X\setminus Y$ in the $q^{-1}$-closure of $Y$ with a totally bounded $q^{-1}$-neighbourghood, then the subspace $Y$ is not $QH$-singular.
\begin{proof}
Choose $\epsilon>0$ such that the $q^{-1}$-ball with center $x$ and radius $\epsilon$ is totally bounded and does not cover the entire subset $Y$. Define $G_n\subseteq Y$ as the subset $\{y\in Y\, |\, q(y, x)> n^{-1}\epsilon\}$. It is clear that this is a $\ll$-increasing sequence of non-empty sets and, since $x$ is an element of $X\setminus Y$ that is contained in the $q^{-1}$-
closure of $Y$, it is a cover of $Y$ that does not contain the set $Y$ itself. It follows from corollary \ref{countablesingular} that $Y$ is not $QH$-singular. 
\end{proof}
\end{prop}

\begin{cor}
Let $(X, q)$ be a uniformly connected, totally bounded quasi-metric space. The space $(X, q)$ is $QH$-singular iff $(X, q^s)$ is compact.
\begin{proof}
That this condition is sufficient was established in corollary \ref{symmcompact}. Now assume that $(X, q)$ is totally bounded and $QH$-singular, but that $(X, q^s)$ fails to be compact. Since $(X, q^s)$ is a totally bounded metric space we know that $(X, q^s)$ is not complete. Let $(\tilde{X}, \tilde{q})$ be the bicompletion of $(X, q)$ and take an $x_0\in \tilde{X}\setminus X$. It is well known (see \cite{Brummer}) that the symmetrisation of the bicompletion of a quasi-uniform space is equal to the completion of the symmetrisation of that space. Because the completion of a totally bounded uniform space is again totally bounded this implies that $(\tilde{X}, \tilde{q})$ is totally bounded. Hence we have that $\tilde{X}$ is a totally bounded $q^{-1}$-neighbourhood of $x_0$. From the fact that $X$ is $\tilde{q}^s$-dense in $\tilde{X}$ we obtain that $x_0$ is an element of the $\tilde{q}^{-1}$-closure of $X$. It follows from proposition \ref{dense} that $(X, q)$ cannot be $QH$-singular and thus we know that $(X, q^s)$ must be compact.
\end{proof}
\end{cor}

\end{document}